\documentclass[J,12pt]{amsart}
\usepackage{epsfig,epsf,amsmath,amssymb}
\newtheorem{theorem}{Theorem}[section]

\newtheorem{definition}{Definition}
\newtheorem{example}[theorem]{Example}
\newtheorem{proposition}[theorem]{Proposition}
\newtheorem{remark}[theorem]{Remark}
\newtheorem{lemma}[theorem]{Lemma}
\begin{document}
\title[]{Essential concepts of digital topology\\
(digital $k$-covering spaces and pseudo $k$-covering spaces)}

\author[]{Sang-Eon Han}
\address[]{Department of Mathematics Education, Institute of Pure and Applied Mathematics\\
	Jeonbuk National University, Jeonju-City Jeonbuk, 54896, Republic of Korea\\
	e-mail address:sehan@jbnu.ac.kr, Tel: 82-63-270-4449.}
\thanks {AMS Classification:54B20, 68U05\\
	Keywords: digital topology, digital image, covering map, unique path lifting property, homotopy lifting theorem, pseudo-covering space, digital covering space, $WL$-$(k_1,k_2)$-isomorphism, local (or $L$-) $(k_1,k_2)$-isomorphism, digital-topological ($DT$-, for brevity) $(k_1,k_2)$-embedding.}

\begin{abstract}
		 	 The present paper focuses on the notions of  covering spaces, pseudo-covering spaces, and their equivalences.
	We discuss something incorrectly mentioned in Boxer's papers and correct them. Indeed, Sections 4-6 (or 4-6) of \cite{B3} are redundant because they have some incorrect assertions on $(k_1,k_2)$-covering spaces or pseudo- $(k_1,k_2)$-covering spaces due to his misunderstanding on Han's papers \cite{H14,H16}.
	In addition, many things in \cite{B3} are duplicated with some results in \cite{H18}.
	In addition, since the papers \cite{P1,P2} also have some defects, we correct and improve them.
	
\end{abstract}

%\date{2013. 3. 16}

\maketitle
\newpage

\section{\bf Terminology for studying digital images}\label{s1}

The present paper will often follow the notation.
${\mathbb N}$ indicates the set of positive integers, ${\mathbb Z}$ means the set of  integers, $\lq\lq$$:=$" is used for introducing a new term, and $\lq\lq$$\subset$" stands for a subset as usual. 
Besides,  for $a, b \in {\mathbb Z}$, 
$[a, b]_{\mathbb Z}:=\{t \in {\mathbb Z}\,\vert \,a \leq t \leq b\}$ and 
$[a, \infty)_{\mathbb Z}:=\{t \in {\mathbb Z}\,\vert \,a \leq t\}$.
We also abbreviate the term $\lq\lq$digital topological" as $DT$.
In addition, $X^\sharp$  indicates the cardinality of the given set $X$.\\

Before 2004 (see \cite{H1,H2,H3}), focusing on the $3$-dimensional real world, all papers in the literature on digital topology dealt with only low dimensional digital image $(X, k)$, i.e., $X \subset {\mathbb Z}^3$ with one of the $k$-adjacency of
${\mathbb Z}^3$ \cite{B2,KR1,R1,R2}. 
For instance, the $2$-adjacency for ${\mathbb Z}$; the $4$- and $8$-adjacency for ${\mathbb Z}^2$; and the $6$-, $18$-, and $26$-adjacency for ${\mathbb Z}^3$ are considered \cite{R1,R2}.\\
After generalizing these $k$-connectivities into those for $n$-dimensional cases, $n \in {\mathbb N}$, since 2004 (\cite{H1,H2,H3}), 
a digital image $(X, k)$ has often been assumed to be a set $X$ in ${\mathbb Z}^n$ with one of the $k$-adjacency relations of
${\mathbb Z}^n$  according to (2.2) below. 
Indeed, the $k$-adjacency relations for $X \subset {\mathbb Z}^n, n \in {\mathbb N}$, were initially established in \cite{H3} (see also \cite{H1,H2,H6,H11}).
 In detail, for a natural number $t$, $1 \leq t \leq n$, the distinct points
$p = (p_i)_{i \in [1,n]_{\mathbb Z}}$  and $q=(q_i)_{i \in [1,n]_{\mathbb Z}}\in {\mathbb Z}^n$
are $k(t, n)$-adjacent (see the page 74 of \cite{H3}) if 
$$\text{at most}\,\,t\,\,\text{of their coordinates  differ by}\,\,\pm1\,\,\text{and the others~coincide.}\eqno(2.1)$$

According to the criterion of (2.1), the digital $k(t, n)$-adjacencies of ${\mathbb Z}^n, n \in {\mathbb N}$, are formulated \cite{H3} (see also \cite{H11}) as follows:
$$k:=k(t,n)=\sum_{i=1}^{t} 2^{i}C_{i}^{n}, \text{where}\,\, C_i ^n:= {n!\over (n-i)!\ i!}. \eqno(2.2)$$
For instance, in ${\mathbb Z}^4$, 
$8$, $32$, $64$, and $80$; in ${\mathbb Z}^5$, 
$10$, $50$, $130$, $210$, and $242$; and 
in ${\mathbb Z}^6$, 
$k(1,6)=12$, $k(2,6)=72$, $k(3,6)=232$, $k(4,6)=472$, $k(5,6)=664$, $k(6,6)=728$-adjacency
 are considered \cite{H3,H11}.\\
Besides, the notation $\lq\lq$$k(t, n)$" of (2.2) indeed gives us much information such as the dimension $n$ and the number $t$ associated with the statement of (2.1).
For instance, we may consider an $8$-adjacency in both ${\mathbb Z}^2$ and ${\mathbb Z}^4$.
Then, in the case that we follow the notation $k:=k(t,n)$, we can easily make a distinction between them by using  $k(2,2)=8$ and $k(1,4)=8$ (see (2.1) and (2.2) and \cite{H17}). Hence  Han's approach using $k:=k(t,n)$ can be helpful to study  more complicated cases such as digital products with a $C$-compatible, a normal, and a pseudo-normal $k$-adjacency \cite{H3,H7,H10,KHL1} and their applications.\\
Hereinafter, a relation set $(X, k)$ is assumed in ${\mathbb Z}^n$ with a $k$-adjacency of (2.2), called a digital image. Namely, we now abbreviate $\lq\lq$digital image $(X, k)$" as $\lq\lq$$(X, k)$" if there is no confusion.\\

Let us now recall a $\lq\lq$digital $k$-neighborhood of the point $x_0$" in $(X, k)$.
For $(X, k)$ and a point $x_0 \in X$, we have often used the notation
$$ N_k(x_0, 1):=\{x \in X \,\vert\,\,x\,\,\text{is}\,\,k\text{-adjacent to}\,\,x_0\} \cup \{x_0\}, \eqno(2.3)$$
which is called a digital $k$-neighborhood of the point $x_0$ in $(X, k)$.
The notion of (2.3) indeed facilitates some studies of digital images, e.g., digital covering spaces, digital homotopy, digital-topological ($DT$-, for brevity) $k$-group \cite{H15,H16}, and so on.
Indeed, the notation of (2.3) comes from the $\lq\lq$digital $k$-neighborhood of $x_0 \in X$ with
 radius $\varepsilon \in {\mathbb N}$" in \cite{H1,H2,H3} (see (2.4) below). 
 More precisely, for $(X,k), X \subset {\mathbb Z}^n$ and two distinct points $x$ and $x^\prime$ in $X$, a {\it simple $k$-path} from $x$ to $x^\prime$ with $l$ elements is said to be a sequence $(x_i)_{i \in [0, l]_{\mathbb Z}}$ in $X$ such that $x_0=x$ and $x_l=x^\prime$ and further, $x_i$ and $x_j$ are
$k$-adjacent if and only if $\vert i-j\vert=1$ \cite{KR1}. We say that a {\it length} of the simple $k$-path, denoted by
$l_k(x, y)$, is the number $l$ \cite{H3}.

\begin{definition} \cite{H1} (see also \cite{H3}) For a digital image $(X, k)$ on ${\mathbb Z}^n$,
	the digital $k$-neighborhood of $x_0 \in X$ with
	radius $\varepsilon$ is defined in $X$ to be the following subset of $X$
	$$N_k(x_0, \varepsilon) = \{x \in X \,\vert\, l_k(x_0, x) \leq
	\varepsilon\}\cup\{x_0\}, \eqno(2.4) $$
	where $l_k(x_0, x)$ is the length of a shortest simple $k$-path from $x_0$ to $x$ and  $\varepsilon\in {\mathbb N}$.
\end{definition}

Owing to Definition 1, it is clear that $N_k(x_0, 1)$ of (2.3) is a special case of $N_k(x_0, \varepsilon)$ of (2.4).
Note that the neighborhood $N_k(x_0, 2)$ plays an important role in establishing 
the notions of a radius $2$-local $(k_1,k_2)$-isomorphism and a radius $2$-$(k_1,k_2)$-covering map which are strongly associated with the
homotopy lifting theorem \cite{H2}  and are essential for calculating digital fundamental groups of digital images $(X,k)$ \cite{BK1,H3}. \\

We say that $(Y, k)$ is $k$-connected \cite{KR1} if for arbitrary distinct points $x, y \in Y$ there is a finite sequence $(y_i)_{i \in [0, l]_{\mathbb Z}}$ in $Y$ such that
$y_0=x$, $y_l=y$ and $y_i$ and $y_j$ are
$k$-adjacent if $\vert\, i-j \,\vert=1$ \cite{KR1}.
Besides, a singleton is assumed to be $k$-connected.\\

As $DT$-versions of the adjacencies of graph products established by Berge \cite{B1} and Harray \cite{H19}, a  normal $k$-adjacency  and $C$-compatible $k$-adjacency of a digital product established by Han \cite{H3,H7,H10} (in detail, see \cite{H17}).
Note that a $C$-compatible and normal $k$-adjacency of a digital product  \cite{H3,H7,H10} are different from the adjacencies of graph products established by Berge \cite{B1} and Harray \cite{H19} (see Remarks 3.4 and 3.6, and Example 3.5 of \cite{H17}).\\

A simple closed $k$-curve (or $k$-cycle)
with $l$ elements in ${\mathbb Z}^n, n \geq 2$, denoted by $SC_k^{n,l}$ \cite{H3,KR1}, $l(\geq 4) \in {\mathbb N}$, is defined as the sequence $(y_i)_{i \in [0, l-1]_{\mathbb Z}}\subset {\mathbb Z}^n$ such that
$y_i$ and $y_j$ are $k$-adjacent if and only if  $\vert\, i-j \,\vert=\pm1(mod\,l)$.
Indeed, the number $l(\geq 4)$ is very important in studying $k$-connected digital images 
from the viewpoint of digital homotopy theory.
The number $l$ of $SC_k^{n,l}$ can be an even or an odd number. Namely, the number $l$ depends on the dimension $n$ and the $k$-adjacencies, as follows:

\begin{remark}  (see (5) of \cite{H13}) The number $l$ of $SC_k^{n,l}$ depends on the dimension $n\in {\mathbb N}\setminus \{1\}$ and the number $l \in \{2a, 2a+1\}, a \in {\mathbb N} \setminus \{1\}$.
\end{remark}

\section{\bf Digital $k$-continuity, digital $k$-isomorphism, local $k$-isomorphism, and radius $2$-local $k$-isomorphism}\label{s4}

Rosenfeld \cite{R1} initially introduced the notion of digital $k$-continuity.
More precisely, Rosenfeld \cite{R2} defined that 
a map $f: (X, k_1) \to (X, k_2)$ is a $(k_1, k_2)$-continuous map if  the image $(f(A), k_2)$ is $k_2$-connected for each  $k_1$-connected subset $A$ of $(X, k_1)$, where $X \subset {\mathbb Z}^3$.
Also, Rosenfeld represented it as follows:
 A function $f: (X,k_1) \to (Y, k_2)$ is $(k_1, k_2)$-continuous if and only if  two points $x$ and $x^\prime$ which are $k_1$-adjacent in $(X,k_1)$ implies that either $f(x)=f(x^\prime)$ or $f(x)$ is $k_2$-adjacent to $f(x^\prime)$, where $X$ and $Y$ are subsets of ${\mathbb Z}^3$.
 Based on this approach, Han \cite{H3} represented the $(k_1, k_2)$-continuity of a map $f: (X, k_1) \to (Y, k_2)$ with the following statement which has been used to study digital covering spaces and pseudo-convering spaces, $DT$-$k$-group theory \cite{H15,H16}, and so forth, where $X \subset {\mathbb Z}^{n_1}$ and $Y \subset {\mathbb Z}^{n_2}$.

\begin{proposition} \cite{H1,H3,H5}
	A~function $g:(X, k_1) \to (Y, k_2)$ is
	$(k_1, k_2)$-continuous if and only if for each point  $x\in X$,
	$g(N_{k_1}(x,1))\subset  N_{k_2}(g(x),1)$, where $X \subset {\mathbb Z}^{n_1}$ and $Y \subset {\mathbb Z}^{n_2}$.
\end{proposition}

Using  $k$-continuous maps,
we establish a {\it digital topological category}, denoted by {\it DTC},
consisting of the following two pieces of data \cite{H5}.\\
$\bullet$ The set of $(X, k)$ as  objects,\\
$\bullet$ For every ordered pair of objects  $(X, k_1)$ and $(Y, k_2)$,
the set of all $(k_1, k_2)$-continuous maps $f:(X, k_1) \to (Y, k_2)$ as morphisms.\\

Based on the digital continuity, since 2005 \cite{H4},
we have often used a {\it $(k_1, k_2)$-isomorphism} for a map between two sets $X \subset {\mathbb Z}^{n_1}$ and $Y \subset {\mathbb Z}^{n_2}$ in \cite{H4,HP1} instead of a {\it $(k_1, k_2)$-homeomorphism} for a map between two sets $X \subset {\mathbb Z}^3$ and $Y \subset {\mathbb Z}^3$ in \cite{B2}, as follows:

\begin{definition} \cite{B2} (see also \cite{H4,HP1})
	Assume $(X, k_1)$ and $(Y, k_2)$ on ${\mathbb Z}^{n_1}$ and ${\mathbb Z}^{n_2}$, respectively. We say that
	a bijection $h: (X, k_1) \to (Y, k_2)$ is a
	  $(k_1, k_2)$-isomorphism if $h$
	is a $(k_1, k_2)$-continuous and the inverse map of $h$,
	 $h^{-1}: Y \to X$, is $(k_2, k_1)$-continuous.
	 Then we denote this isomorphism by 
	 $(X, k_1) \approx_{(k_1, k_2)}(Y, k_2)$.
	If $n_1 = n_2$ and $k_1 = k_2$, then we call it a {\it $k_1$-isomorphism}.
\end{definition}

Compared with a $(k_1, k_2)$-isomorphism,
a local $(k_1, k_2)$-isomorphism which is weaker than a  $(k_1, k_2)$-isomorphism was defined, as follows. 
The notion plays an important role in establishing a $(k_1, k_2)$-covering space \cite{H2,H3} and a pseudo-$(k_1, k_2)$-covering space \cite{H9,H14}, $DT$-$(k_1, k_2)$-embedding theorem \cite{H13}, and so forth.

\begin{definition} \cite{H2,H5}
	For two digital images $(X, k_1)$ on ${\mathbb Z}^{n_1}$ and $(Y, k_2)$ on ${\mathbb Z}^{n_2}$, consider a  map  $h: (X,k_1) \to  (Y,k_2)$. Then~the map $h$ is said to be a local (or $L$-, for brevity) $(k_1, k_2)$-isomorphism 
	if  for every $x \in X$, $h$ maps $N_{k_1}(x, 1)$ $(k_1, k_2)$-isomorphically onto $N_{k_2}(h(x), 1)$, i.e., the restriction map $h\vert_{N_{k_1}(x, 1)}:
	N_{k_1}(x, 1) \to N_{k_2}(h(x), 1)$ is a $(k_1, k_2)$-isomorphism.
	If $n_1 = n_2$ and $k_1 = k_2 $, then the map $h$ is called a local (or $L$-) $k_1$-isomorphism.
\end{definition}

In view of Definition 3, a local $(k_1, k_2)$-isomorphism is obviously a $(k_1, k_2)$-continuous map. However, the converse does not hold.
Compared with a local $(k_1, k_2)$-isomorphism, the following is defined.

\begin{definition} \cite{H2} (Radius $2$-local $(k_1,k_2)$-isomorphism)
	Assume a map  $h:(X, k_1) \to (Y, k_2)$, where $X \subset {\mathbb Z}^{n_1}$ and $Y \subset {\mathbb Z}^{n_2}$.
	Then we say that it is a radius $2$-local $(k_1,k_2)$-isomorphism if 
	for each point $x\in X$ the restriction $h$ to $N_{k_1}(x, 2)$, i.e., 
	$h\vert_{N_{k_1}(x, 2)}:N_{k_1}(x, 2) \to N_{k_2}(h(x), 2)$ is a 
	$(k_1,k_2)$-isomorphism.
	In particular, in the case of $k:=k_1=k_2$, we call it a radius $2$-local
	$k$-isomorphism.
	\end{definition}

Comparing Definitions 3 and 4, we obtain the following:
\begin{remark} A radius $2$-local $(k_1,k_2)$-isomorphism is clearly a local $(k_1, k_2)$-isomorphism. However, the converse does not hold.
\end{remark}

\begin{example}  Consider the map $p:({\mathbb Z},2) \to SC_k^{n,l}:=(s_i)_{i \in [0,l-1]_{\mathbb Z}}$ given by $p(t)=s_{t(mod\,l)}$.\\
		(1) In the case
		of $l\geq 5$, the map $p$ is a
	 radius $2$-local $(2,k)$-isomorphism.\\	
	(2) In the case of $l=4$, the map $p$ is not a radius $2$-local $(2,k)$-isomorphism but an $L$-$(2,k)$-isomorphism.	
	\end{example}

A radius $2$-local $(k_1,k_2)$-isomorphism  strongly plays an important role in studying digital images from the viewpoint of digital homotopy theory and digital covering theory.
In detail, this notion is essential for establishing the homotopy lifting theorem \cite{H2} which can be used to calculate digital fundamental groups \cite{BK1,H2,H3}.\\
\begin{remark} (Utilities of  radius $2$-local $(k_1,k_2)$-isomorphism)\\
	$\bullet$ Establishment of the homotopy lifting theorem \cite{H2}\\
	$\bullet$ Calculation of digital fundamental groups \cite{BK1,H2,H3}
	\end{remark}

\section{\bf Remarks on digital $k$-covering spaces and  pseudo-$k$-covering spaces: Unique path lifting property and digital homotopy lifting theorem}\label{s6}

In 2004 and 2005, Han \cite{H1,H2,H3} initially introduced the notion of a digital $(k_1, k_2)$-covering space and investigated various properties of it, 
which has facilitated the study of digital topology. 
In 2012, to make a classical digital covering space broader,  the first version of a pseudo-$(k_1, k_2)$-covering map was introduced in \cite{H9}. Since the notions of a digital $(k_1,k_2)$-covering map (see Definition 5 below) and a pseudo-$(k_1,k_2)$-covering  (see Definition 8 below) were designed to be different from each other (compare the condition (3) of Definition 5 and that of Definition 8).
However, owing to the conditions (1) of Definitions 5 and 8, 
it recently turns out that these two notions are equivalent \cite{P1} (see also \cite{H18}).
To confirm this equivalence, we first recall these two notions (see Definitions 5 and 8). Eventually, the condition (1) for the original version of a pseudo-covering space of Definition 8 \cite{H9}  was revised (compare the conditions (1) of Definitions 5 and 10 in the present paper).
Hence  we now note that the pseudo-$k$-covering map of Definition 10 is broader than the digital $k$-covering map  of Definition 5 (see \cite{H18}). \\

Based on the digital $(k_1,k_2)$-covering map in \cite{H1,H2,H3}, if we take $\varepsilon=1=\delta$, we obtain the following.

\begin{definition} \cite{H5} (see Definition 5 and Remark 1 on page 221 of \cite{H5})
	Let $(E, k_1)$ and $(B, k_2)$ be digital images on  ${\mathbb Z}^{n_1}$ and  ${\mathbb Z}^{n_2}$, respectively.
	Let $p: E \to B $ be a surjection such that for any $b \in B$,\\
	(1) for some index set $M$, $p^{-1}(N_{k_2}(b, 1)) = \bigcup\limits_{i \in M}N_{k_2}(e_i, 1)$ with $e_i \in p^{-1}(b):=p^{-1}(\{b\})$;\\
	(2) if  $i, j \in M$ and $i \neq j$, then
	$N_{k_2}(e_i, 1)\cap N_{k_2}(e_j, 1)$ is an empty set; and\\
	(3) the restriction of $p$ to $N_{k_2}(e_i, 1)$ from $N_{k_1}(e_i, 1)$ to $N_{k_2}(b, 1)$ is a $(k_1, k_2)$-isomorphism for all $i \in M$.\\
	Then the map $p$ is called a digital $(k_1, k_2)$-covering map, $(E, p,
	B)$ is said to be a digital $(k_1, k_2)$-covering and $(E, k_1)$ is called a digital $(k_1, k_2)$-covering space over $(B, k_2)$.
\end{definition}   
                                                                
Hereinafter, we abbreviate $\lq\lq$digital $(k_1, k_2)$-covering map" as $(k_1, k_2)$-covering map.
Note that the  $(k_1, k_2)$-covering map of Definition 5 is usually called a radius $1$-$(k_1, k_2)$-covering map \cite{H5} which can be assumed to be a typical $(k_1, k_2)$-covering map. 
Definition 5 is a special case of the following:

\begin{definition} \cite{H1,H2,H3,H5}
	Let $(E, k_1)$ and $(B, k_2)$ be digital images in ${\mathbb Z}^{n_1}$
	and ${\mathbb Z}^{n_2}$, respectively
	and let $p: E \to B $ be a $(k_1, k_2)$-continuous surjection.
	Suppose that for every $b \in B$ there exists $\varepsilon \in {\mathbb N}$ such that
		$(DC 1)$ for some index set $M$,	
	$p^{-1}(N_{k_2}(b, \varepsilon)) = \cup _{i \in M}N_{k_1}(e_i, \varepsilon)$ with $ e_i \in p^{-1}(\{b\})$;\\	
	$(DC 2)$ if  $i, j \in M$ and $i \neq j$, then
	$N_{k_1}(e_i, \varepsilon)\cap N_{k_1}(e_j, \varepsilon) = \phi$; and\\	
	$(DC 3)$ the restriction map
	$p\vert_{N_{k_1}(e_i, \varepsilon)}:N_{k_1}(e_i, \varepsilon) \to N_{k_2}(b, \varepsilon)$
	is a $(k_1, k_2)$-isomorphism for all $i \in M$.\\	
	Then the map $p$ is called a $(k_1, k_2)$-covering map
	and $(E, p, B)$ is said to be a $(k_1, k_2)$-covering.
\end{definition}

In Definition 6, if $\varepsilon=2$, then we call the map $p$ a  radius $2$-$(k_1, k_2)$-covering map \cite{H5}.

Let us consider some reason why Han \cite{H3,H5} took the number $\varepsilon$ for the notations $N_{k_1}(e_j, \varepsilon)$ and $N_{k_2}(b, \varepsilon)$ of Definition 6.
In Definition 6, $N_{k_2}(b, \varepsilon)$
is called an {\it elementary $k_2$-neighborhood} of $b$ with radius $\varepsilon$.
The collection $\{N_{k_1}(e_i,  \varepsilon) \,\vert\, i \in M \}$ is
a partition of $p^{-1}(N_{k_2}(b, \varepsilon))$ into slices.\\
In Definition 6, if we take $\varepsilon=1$, then we obtain Definition 5. Namely, Definition 5 is a special case of Definition 6 and it is a typical $(k_1, k_2)$-covering map (see Definition 5 and Remark 1 on page 221 of \cite{H5}).

Let us now mention that Definition 6 is also useful. Indeed, in Definition 6 we are only interested in the case $\varepsilon \in \{1,2\}$.
In the case of $\varepsilon=1$, we have  Definition 5.
One of the important things is that in the case of
$\varepsilon=2$, the $(k_1, k_2)$-covering map of Definition 6 also plays an important role in establishing the digital homotopy lifting theorem which is essential for calculating digital fundamental groups of digital images.
Indeed, using both the radius $2$-local $(k_1, k_2)$-isomorphism of Definition 4 and 
the $(k_1, k_2)$-covering map of Definition 5 we can observe the following:

 \begin{remark} A $(k_1, k_2)$-covering map of Definition 5 with a radius $2$-local $(k_1, k_2)$-isomorphism is equivalent to the $(k_1, k_2)$-covering map of Definition 6 with $\varepsilon=2$.
 	Indeed, we often abbreviate $\lq\lq$$(k_1, k_2)$-covering map of Definition 6 with $\varepsilon=2$" as $\lq\lq$radius $2$-$(k_1, k_2)$-covering map (see \cite{H2}).
 	 \end{remark}

\begin{remark} In the paper \cite{BK1}, Boxer et al. stated that a $(k_1, k_2)$-covering map of Definition 6 is equivalent to the $(k_1, k_2)$-covering map of Definition 5 as a proposition (see Proposition 2.4 in \cite{BK1}). 
However, it is not a proposition made in \cite{BK1} but just a definition as stated in Definition 5.
	\end{remark}

As an equivalent version of a $(k_1, k_2)$-covering map of Definition 5, the following is obtained, which can play an important role in studying digital covering spaces.
\begin{proposition} (1) An $L$-$(k_1, k_2)$-isomorphic surjection is equivalent to a $(k_1, k_2)$-covering map \cite{PZ1}.\\
	(2) For $k_1$-connected digital image $(E,k_1)$ and $k_2$-connected digital image $(B,k_2)$, an $L$-$(k_1, k_2)$-isomorphism $p:(E,k_1) \to (B,k_2)$ is equivalent to a $(k_1, k_2)$-covering map of $p$ \cite{H14}.	
\end{proposition}

Since both the unique digital lifting theorem and the digital homotopy lifting theorem have been used to study digital covering spaces, let us now review them and related results, as follows.\\
For three digital images $(E, k_1)$ on ${\mathbb Z}^{n_1}$, $(B, k_2)$ on ${\mathbb Z}^{n_2}$,
and $(X, k)$ on ${\mathbb Z}^{n}$,
let $p: E \to B $ be a  $(k_1, k_2)$-continuous map.
For a $(k, k_2)$-continuous map $f:(X, k) \to (B, k_2)$,
as the digital analogue of the lifting in algebraic topology,
we say that a {\it digital lifting} of $f$ is a $(k, k_1)$-continuous map $\tilde f :X  \to E$
such that $p \circ \tilde f = f$ \cite {H3}.
We now recall the {\it unique digital lifting theorem} in \cite{H3}, as follows.

\begin{lemma} \cite{H3} For pointed digital images $((E, e_0), k_1)$ on ${\mathbb Z}^{n_1}$
	and $((B, b_0), k_2)$ on ${\mathbb Z}^{n_2}$,
	let $p: (E, e_0) \to (B, b_0)$ be a pointed $(k_1, k_2)$-covering map.
	Any $k_2$-path $f: [0, m]_{\mathbb Z} \to B$ beginning at $b_0$
	has a unique digital lifting to a $k_1$-path $\tilde f$ in $E$ beginning at $e_0$.
\end{lemma}

Moreover, the following {\it digital homotopy lifting theorem} was introduced in \cite{H2},
which plays an important role in studying digital covering theory.

\begin{lemma} \cite{H2} Let $((E, e_0), k_1)$ and $((B, b_0), k_2)$ be pointed digital images.
	Let $p: (E, e_0) \to (B, b_0)$ be a radius $2$-$(k_1, k_2)$-covering map.
	For $k_1$-paths
	$g_0, g_1$ in $(E, e_0)$ that start at $e_0$, if there is  a $k_2$-homotopy in $B$
	from $p \circ g_0 $ to $ p \circ g_1$ that holds the endpoints fixed,
	then  $g_0$ and $ g_1$ have the same terminal point,
	and there is a $k_1$-homotopy in $E$ from $g_0$ to $g_1$ that holds the endpoints fixed.
\end{lemma}

The paper \cite{H9} defined the following notion which is weaker than an $L$-$(k_1, k_2)$-isomorphism.

\begin{definition} \cite{H9} For two digital images $(X, k_1)$ on ${\mathbb Z}^{n_1}$ and
	$(Y, k_2)$ on ${\mathbb Z}^{n_2}$, a map $h: X \to Y$ is called a
	weakly local (WL-, for brevity) $(k_1, k_2)$-isomorphism
	if  for every $x \in X$, $h$ maps $N_{k_2}(x, 1)$ $(k_1, k_2)$-isomorphically onto $h(N_{k_1}(x, 1))\subset (Y, k_2)$, i.e., the restriction map $h\vert_{N_{k_1}(x, 1)}:
	N_{k_1}(x, 1) \to h(N_{k_1}(x, 1))$ is a $(k_1, k_2)$-isomorphism.
	In particular, if $n_1 = n_2$ and $k:=k_1 = k_2$, then the map $h$ is called a weakly local $k$-isomorphism (or a $WL$-$k$-isomorphism).
\end{definition}

Using this notion, the paper \cite{H9} defined  the notion of a pseudo-$(k_1, k_2)$-covering space, as follows:
\begin{definition} \cite{H9} (The first version of a pseudo-$(k_1, k_2)$-covering map)
	Let $(E, k_1)$ and $(B, k_2)$ be digital images on  ${\mathbb Z}^{n_1}$ and  ${\mathbb Z}^{n_2}$, respectively.
	Let $p: E \to B $ be a surjection such that for any $b \in B$,\\
	(1)-(2): the conditions (1) and (2) of Definition 5;  and\\
	(3) the restriction of $p$ to $N_{k_1}(e_i, 1)$ from $N_{k_1}(e_i, 1)$ to $N_{k_2}(b, 1)$ is a $WL$-$(k_1, k_2)$-isomorphism for all $i \in M$.\\
	Then the map $p$ is called a pseudo-$(k_1, k_2)$-covering map, $(E, p,
	B)$ is said to be a pseudo-$(k_1, k_2)$-covering and $(E, k_1)$ is called a pseudo-$(k_1, k_2)$-covering space over $(B, k_2)$.
	\end{definition}

Based on the notion of a pseudo-$(k_1, k_2)$-covering space, 
the paper \cite{H9} suggested an example for a pseudo-$(2,k)$-covering map with the map $$p:({\mathbb Z}^+=[0, \infty)_{\mathbb Z}, 2) \to SC_k^{n,l}:=(s_i)_{i \in [0, l-1]_{\mathbb Z}} \eqno (3.1), $$
defined by $p(t)=s_{t(mod\,l)}$.\\
Indeed, the paper \cite{H9} made a mistake to take this map as a pseudo-$(2,k)$-covering map (see Proposition 3.6 below and Remark 3.5 of \cite{H18}).
By contrast to the condition (1) of Definition 8, the map $p$ of (3.1) is not a pseudo-$(2,k)$-covering map, as follows:

\begin{proposition} (Proposition 3.2 of \cite{P1})
	The map $p:({\mathbb Z}^+,2) \to SC_k^{n,l}:=(c_i)_{i \in [0, l-1]_{\mathbb Z}}, l \geq 4$,  is not a pseudo-$(2,k)$-covering map of Definition 8.
\end{proposition}

However, the proof of Proposition 3.6 (or Proposition 3.2 of \cite{P1}) is clearly incorrect.
Thus the paper \cite{H18} corrected the defect (see Remark 3.5 of \cite{H18}), as follows:

\begin{remark} (Correction of the proof of Proposition 3.4 in the present paper or Proposition 3.2 of \cite{P1} and the lines from 13-14 on page 51 of \cite{P2})\\
		Based on the proof of Proposition 3.2 of \cite{P1} and some part in \cite{P2}, we can clearly see some errors in the assertions of \cite{P1,P2}.
		First of all, let us see the proof of Proposition 3.2 of \cite{P1} which is incorrect.
		In detail, see the part with lines 9-10 from the bottom of page 214 of \cite{P1} (see  Remark 3.5 of \cite{H18}).\\
			In addition, a similar mistake appears in the paper \cite{P2}.
	See the part from the line 13-14 which is incorrect.
	Indeed, the part should be written $6k-1$ instead of $5k$ (in detail, see Remark 3.5 of \cite{H18}).
	\end{remark}

Owing to the condition (1) of Definition 8, the following is obtained.

\begin{theorem} (See Theorem 3.3 of \cite{P1}) A pseudo-$(k_1, k_2)$-covering map Definition 8 is a $(k_1, k_2)$-covering map of Defintion 5.
\end{theorem}

\begin{remark}
The proof of Theorem 3.8  was also improved in \cite{H18} (see the proof of Theorem 3.7 of \cite{H18}).
\end{remark}

\section{\bf Remarks on Boxer's misconception on Han's papers \cite{H9,H14} and a revision of the condition (1) of Definition 8 for a pseudo-$(k_1, k_2)$-covering space}\label{s6}

In  Section 4 (or $4$) of Boxer's manuscript \cite{B3}, we can clearly observe Boxer's recklessness. More precisely, in the part, Boxer incorrectly referred to Han's result in  the part (3) of Theorem 3.15 of \cite{H14}.
To make the paper self-contained, let us recall the notion of a $PL$-$(k_1, k_2)$-isomorphism in \cite{H1,H14}, as follows:

\begin{definition} (\cite{H1,H14}) For two digital images $(X, k_1)$ on ${\mathbb Z}^{n_1}$ and $(Y,k_2)$ on ${\mathbb Z}^{n_2}$, a $(k_1, k_2)$-continuous map $h: X \to
	Y$ is called a pseudo-local ($PL$-, for brevity) $(k_1, k_2)$-isomorphism if
	for every point $x \in X$, $h(N_{k_1}(x, 1))$ is $k_2$-isomorphic with $N_{k_2}(h(x), 1)$.
	If $n_1 = n_2$ and $k_1 = k_2 $, then the map $h$ is called a $PL$-$k_1$-isomorphism.
\end{definition}

 Consider a map $f:(X,k_1) \to (Y, k_2)$. Then we have the following result (see Theorem 3.15(3) of \cite{H14}).\\

$$ 
\left \{\aligned 
&\text{Neither of a} \,\,PL\text{-}(k_1, k_2)\text{-isomorphism of the map}\,\,f\\
&\text{and a}\,\, WL\text{-}(k_1, k_2)\text{-isomorphism of the map}\,\,f\\
& \text{implies the other}.
\endaligned
\right\} \eqno (4.1)
$$

\begin{remark} (Correction of Boxer's assertion in \cite{B3})
As for the result of (4.1), Boxer agreed that the assertion of (4.1) is correct. However, Boxer \cite{B3} wrote that the proof for (4.1) as in \cite{H14} can be insufficient  due to an absence of the $k_2$-connectedness of $(Y, k_2)$.
Indeed, this observation is due to Boxer's recklessness because the paper \cite{H14} wrote the assertion with the condition that the digital image 
$(Y, k_2)$ is $k_2$-connected.
To confirm the fact, we can see the last line of Abstract of the paper \cite{H14}.
Han \cite{H14} asserts all results in \cite{H14} including (4.1) with the assumption that each digital image $(X,k)$ in the paper \cite{H14} is $k$-connected unless otherwise stated.
\end{remark}

As for Boxer's misconception of \cite{P1}, one of the important things is the following:

\begin{remark} 
	Even though the proof of Proposition 3.2 of \cite{P1} is clearly incorrect as mentioned in Remark 3.7, Boxer \cite{B3} missed this defect in the proof of Proposition 3.2 of \cite{P1} and approved this approach of \cite{P1}. 
	To confirm this fact, see Remark 3.5 of \cite{H18}. Indeed, Han \cite{H18} also improved the proof of  Proposition 3.2 of \cite{P1}.
\end{remark}

Motivated by Proposition 3.6 and Theorem 3.8, let us now revise the condition (1) of Definition 8 to make a distinction from a $(k_1, k_2)$-covering map. 

\begin{definition}\cite{H18} (Revision of Definition 8 for a pseudo-$(k_1, k_2)$-covering space) 	Let $(E, k_1)$ and $(B, k_2)$ be digital images in  ${\mathbb Z}^{n_0}$  and  ${\mathbb Z}^{n_1}$, respectively.
	Let $p: E \to B $ be a surjection such that for any $b \in B$,\\
	(1) for some index set $M$, $\bigcup\limits_{i \in M}N_{k_1}(e_i, 1)$ is a subset of $p^{-1}(N_{k_2}(b, 1))$  with $e_i \in p^{-1}(\{b\})$
	and the  other two conditions are equal to the conditions (2) and (3) of Definition 8.\\
	Then the map $p$ is called a  pseudo-$(k_1, k_2)$-covering map and $(E, k_1)$ is called a  pseudo-$(k_1, k_2)$-covering space over $(B, k_2)$.
	In particular, in the case of $n_0=n_1$ and $k:=k_1=k_2$, we use the term 
	$\lq$pseudo-$k$-covering'  instead of  $\lq\lq$pseudo-$(k_1, k_2)$-covering".
	\end{definition}

Hereinafter, we will follow Definition 10 for a pseudo-$(k_1, k_2)$-covering space as a revised version of Definition 8.
Then we have some utilities of it and many good examples in \cite{H18}.

	Based on the pseudo-$(k_1, k_2)$-covering space of Definition 10, 
the map $p$ of (3.1) is now a pseudo-$(2,k)$-covering map \cite{H18}.

With Definition 10, we note that a pseudo-$(k_1, k_2)$-covering space is now broader than a  $(k_1, k_2)$-covering space \cite{H18}.
For instance, we note that the map $p$ of (3.1) is not a $(2,k)$-covering map but a pseudo-$(2,k)$-covering map based on Definition 10.
 Besides, many examples discriminating between these two notions were referred to in \cite{H18} (see Theorem 4.2 and Example 4.3 of \cite{H18}). 
 Hereinafter, we follow the notion of Definition 10 for a pseudo-$(k_1, k_2)$-covering space.

After comparing the condition (1) of Definition 5 and that of Definition 10,
the following is obtained.  

\begin{theorem} \cite{H18}
 A pseudo-$(k_1, k_2)$-covering map of Definition 10 is broader than a $(k_1, k_2)$-covering map of Definition 5.
	\end{theorem}

\begin{remark}
	The paper \cite{H18} suggests some new relationships between a $WL$-$(k_1, k_2)$-isomorphism and a pseudo-$(k_1, k_2)$-covering map of Definition 10. 
\end{remark}

As mentioned in \cite{H18}, a pseudo-$(k_1, k_2)$-covering map of Definition 10 supports neither a unique path lifting property nor the homotopy lifting theorem of Lemma 3.5.

\section{\bf Remarks on a $DT$-$(k_1, k_2)$-embedding and some incorrect assertions in Boxer's manuscript \cite{B3}}\label{s7}

As for some assertions in section 6 (or $6$) of \cite{B3}, there are some parts which should be corrected.
Since the notion of a $DT$-$(k_1, k_2)$-embedding can play a crucial role in applied topology, a recent paper \cite{H13} proposed it, as follows:

\begin{definition} \cite{H13} ($DT$-$(k_1, k_2)$-embedding) 
	Assume $(X, k_1:=k(t_2,n_1))$, $X \subset {\mathbb Z}^{n_1}$ and $(Y, k_2:=k(t_2,n_2))$,
	$Y \subset {\mathbb Z}^{n_2}$ such that there is
	an arbitrary $(k_1, k_2)$-isomorphism $h:(X, k_1) \to (h(X), k_2)\subset (Y, k_2)$.
	Then
	 the map $h$ is said to be a $DT$-$(k_1,k_2)$-embedding of $(X, k_1)$ into $(Y,k_2)$ or  $(X, k_1)$ is a  $DT$-$(k_1, k_2)$-embedding into $(Y, k_2)$ {\it w.r.t.} the $(k_1, k_2)$-isomorphism $h$.\\ 
	In the case of $k:=k_1=k_2$, we call it a $DT$-$k$-embedding instead of  $\lq\lq$$DT$-$(k_1, k_2)$-embedding".
	\end{definition}

\begin{theorem} (Theorem 3.20 of \cite{H14}) Given two $SC_{k_i}^{n_i, l_i}, i \in \{1,2\}$, 
	$SC_{k_1}^{n_1, l_1}$ is a $DT$-$(k_1, k_2)$-embedding into $SC_{k_2}^{n_2, l_2}$  if and only if $l_1=l_2$.
\end{theorem}

As for Theorem 5.1, Boxer \cite{B3} mentioned some relationships between a $DT$-$(k_1, k_2)$-embedding and a $(k_1, k_2)$-homotopy equivalence, we can see that
 Boxer's observation seems to be too far. Namely, we need to focus on a  $DT$-$(k_1, k_2)$-embedding. One can easily see that 
 neither of a $DT$-$(k_1, k_2)$-embedding and a $(k_1, k_2)$-homotopy equivalence implies the other. 
 
 \begin{remark}
 	The paper \cite{H18} has several results mentioning some relationships among many variants of a  $(k_1, k_2)$-covering map and a pseudo-$(k_1, k_2)$-covering map.
 \end{remark}
 
 In \cite{B3}, Boxer referred to an example for a pseudo-$4$-covering map suggested in Han's paper \cite{H9} (see Example 4.3(4) of \cite{H9} and Remark 6.6 of \cite{B3}).
 To make the paper self-contained, let us recall the map $p$ in \cite{H9}
  Let $C:=\{c_0=(0,0), c_1=(1,0), c_2=(1,1), c_3=(0,1)\}$ and $D:=\{d_0=c_0, d_1=c_1, d_2=c_2, d_3=c_2\}$ 
 and consider both $C$ and $D$ with $4$-adjacency, i.e., $(C,4)$ and $(D,4)$.
 Then consider the map 
 $$p:(C,4) \to (D,4)\,\,\text{given by}\,\,p(c_i)=d_i, i \in [0,3]_{\mathbb Z}.\eqno(5.1)$$
  Then, based on the pseudo-$(k_1,k_2)$-covering of Definition 8,
 Han \cite{H9} wrote that the map $p$ is not a pseudo-$4$-covering map.
 Indeed, Boxer \cite{B3} also agreed that this result is correct (see Remark 6.6 of \cite{B3}).
 However, he added some strange comment that the map $p$ is not $4$-continuous map.
 Hence the map $p$ can be an unsuitable example for disproving a pseudo-$4$-covering map. \\
 As for this comment (or speculation and guess), Boxer's opinion seems to be too far from Han's intention.
 Indeed, the map $p$ of (5.1) is one of the good examples for a non-pseudo-$4$-covering map. We need not refer to the $(4, 8)$ or $4$-continuous map of $p$.
 	Namely, we need to press on the non-pseudo-$4$-covering map of $p$.
 	In detail, to examine if the map $p$ is a pseudo-$4$-covering, we need to prove the following:\\
 	After proving a surjection of $p$, we need only to show that the map $p$ satisfies the conditions (1)-(3) of Definition 8.
 Let us examine if the map $p$ of (5.1) is a  pseudo-$4$-covering map.
 First, it is clear that it is a surjection.
 Second, the map $p$ clearly satisfies the conditions (1) and (2) of Definition 8. 
 Finally, we need to ask if the map $p$ satisfies the condition (3) of  Definition 8.
 In detail, based on Proposition 2.1, we can see that the map $p$ is not a $WL$-$4$-isomorphism (see the point $c_0$), i.e., 
  $$ 
 \left \{\aligned
 & p(N_4(c_0, 1))=p(\{c_0, c_1, c_3\})=\{d_0, d_1, d_3\} \\
 & \nsubseteq N_4(p(c_0), 1)=N_4(d_0, 1)=\{d_0, d_1, d_2\},
 \endaligned
 \right\} \eqno (5.2) 
 $$
  which implies that the map $p$ is not a $WL$-$4$-isomorphism.
  Hence Boxer \cite{B3} should not have considered anything which is not related to a 
 $WL$-$4$-isomorphism. Indeed, the map $p$ is clearly not a $4$-continuous map.
 Of course, the map $p$ is trivially not a $WL$-$(4,8)$-isomorphism because 
 $N_8(d_i, 1)=D$ for each $d_i \in D, i \in [0,3]_{\mathbb Z}$. 
 Besides, even though we follow the new version of a pseudo-$(k_1,k_2)$-covering map of Defintion 10, the map of (5.1) is not a pseudo-$4$-covering map.

 \section{\bf Further remarks}\label{s8}
 
 We have discussed many defects in Boxer's paper \cite{B3}.  We can confirm that 
 there are several maps which are equivalent to a $(k_1,k_2)$-covering map.
 In particular, the paper \cite{H18} corrected and improved all results relating to a covering map, a pseudo-covering map and related maps in the literature on digital covering spaces.\\

 {\bf Conflicts of Interest}: The author declares no conflict of interest.\\

\end{document}